\documentclass[12pt]{article}

\usepackage{amsmath,graphicx}
\usepackage[pdfpagemode=None,colorlinks=true,linkcolor=blue,citecolor=blue]{hyperref}
\usepackage{amsfonts,amssymb}
\usepackage{microtype}

\newtheorem{prop}{Proposition}[section]
\newtheorem{rema}[prop]{Remark}
\newtheorem{defi}[prop]{Definition}
\newtheorem{lemm}[prop]{Lemma}
\newtheorem{theo}{Theorem}
\newtheorem{coro}[prop]{Corollary}
\newtheorem{exem}[prop]{Example}

\newcommand{\C}[1][]{\ensuremath{{\mathbb{C}^{#1}} }}
\newcommand{\R}[1][]{\ensuremath{{\mathbb{R}^{#1}} }}
\renewcommand{\S}[1][]{\ensuremath{{\mathbb{S}^{#1}} }}
\renewcommand{\H}[1][]{\ensuremath{{\mathbb{H}^{#1}} }}
\newcommand{\G}{\mathbb{G}}
\newcommand{\J}{\mathbb{J}}
\newcommand{\Om}{\Omega}
\newcommand{\QED}{\normalsize {}%
\nolinebreak \hfill $\blacksquare$ \medbreak \par}

\newcommand{\<}{\langle}
\renewcommand{\>}{\rangle}
\newcommand{\ga}{\gamma}
\newcommand{\pa}{\partial}

\newcommand{\al}{\alpha}
\newcommand{\eps}{\epsilon}
\newcommand{\te}{\theta}

\newcommand{\la}{\lambda}
\newcommand{\be}{\beta}
\newcommand{\si}{\sigma}

\date{}

\author{Henri Anciaux, Brendan Guilfoyle, Pascal Romon}
\title{Minimal Lagrangian surfaces in the tangent bundle of a Riemannian surface }
\begin{document}
\maketitle

\begin{abstract}
Given an oriented Riemannian surface $(\Sigma, g)$, its tangent bundle
$T\Sigma$ enjoys a natural pseudo-K\"{a}hler structure, that is the combination of a complex structure
$\J$, a pseudo-metric $\G$ with neutral signature and a
symplectic structure $\Om$. We give a local classification
of those surfaces of $T\Sigma$ which are both Lagrangian with respect to $\Om$ and minimal with respect to
$\G$. We first show that if $g$ is non-flat, the only such surfaces are affine normal bundles over geodesics.
In the flat case there is, in contrast, a large set of Lagrangian minimal surfaces, which is described explicitly.
As an application, we show that motions of surfaces in $\R^3$ or $\R^3_1$ induce
Hamiltonian motions of their normal congruences, which are Lagrangian surfaces in $T\S^2$ or $T \H^2$ respectively. 
We relate the area of the congruence to a second-order functional $\mathcal{F}=\int
\sqrt{H^2-K}\, dA$ on the original surface.

\medskip

\end{abstract}

\section*{Introduction}

It has recently been observed (\em cf. \em \cite{GK1},\cite{GK2}) that
the tangent bundle $T \Sigma$ of an oriented Riemannian surface $(\Sigma,g,j)$
with metric $g$ and complex structure $j$ enjoys a rich structure:
besides the symplectic form $\Om$ obtained by pulling back the
canonical symplectic form of $T^*\Sigma,$ it can be endowed with a
natural complex structure $\J$ depending on the complex structure
$j$; next we may define a symmetric 2-tensor $\G$ by combining
$\Omega$ and $\J$ in the formula $\G(.,.)=\Om(\J.,.)$. It turns
out that $\G$ is a pseudo-Riemannian metric on $T\Sigma$ with
signature $(2,2)$ and that the complex structure $\J$ is parallel
with respect to $\G$; in other words we have a pseudo-K\"{a}hler
structure on $T \Sigma$. Of particular interest is the case of $\Sigma$ 
being the
two-sphere $\mathbb{S}^2$, since $T\S^2$ can be naturally identified
with the space of oriented lines of Euclidean three-space $\R^3$.
Moreover, under this identification, a two-parameter family of
lines in $\R^3$ ---thus a surface in $T\S^2$--- is Lagrangian if
and only if the lines are normal to some surface of $\R^3$.

\bigskip

Natural objects of study in K\"ahler geometry are minimal Lagrangian submanifolds 
(\em cf. \em \cite{S},\cite{SW}). In the particular case of K\"{a}hler-Einstein manifolds, 
it is a remarkable fact that the mean curvature vector $\vec{H}$ 
of a Lagrangian submanifold $L$ of dimension $n$ is related 
by the formula $\vec{H} = \frac{1}{n} J \nabla \beta$ to a function $\beta$,
the \emph{Lagrangian angle}, defined on the submanifold.
A striking consequence of this formula is that a
 Lagrangian
submanifold is minimal if and only if it has constant Lagrangian
angle. From the analytical viewpoint this structure reduces the
order of the corresponding PDE from 2 to 1.
 In the context of Calabi-Yau geometry, these
submanifolds are in addition calibrated and thus minimizers, and
are called \em Special Lagrangian \em submanifolds (\emph{cf.}~\cite{HL}).

\bigskip

In this paper we give a local classification of
 minimal Lagrangian
surfaces in $(T \Sigma,\J,\G,\Om)$. It turns out that that the picture is
strongly contrasted between on the one hand, the non-flat case,
which is very rigid in the sense that the only non-trivial minimal
Lagrangian surfaces are the normal bundles over a
geodesic of $\Sigma$ (Theorem \ref{one}), and on the other hand
the flat case, where there exists a variety of minimal Lagrangian
surfaces. Moreover, in Euclidean $4$-space endowed with the
standard pseudo-K\"{a}hler metric of signature $(2,2),$ we can attach
to a Lagrangian surface
 a kind of Lagrangian angle function, still satisfying the formula
 $\vec{H} =\frac{1}{2} \J \nabla \beta$, and thus
 whose constancy characterizes
 minimal Lagrangian surfaces. Finally,
 the underlying partial differential equation is linear and thus
 can be explicitly integrated (Theorem \ref{two}).

\bigskip

Another important class of Lagrangian surfaces are those which are
critical points of the area functional restricted to Hamiltonian
variations (\em cf. \em \cite{SW}). The corresponding Euler equation
is the vanishing of the divergence of the mean curvature vector
(for the induced metric). We give some non-trivial examples
of Hamiltonian stationary Lagrangian surfaces in $(T\Sigma,\J,\G,\Om)$.
 In the special case
of $T\S^2,$ we already know that these Hamiltonian stationary
Lagrangian surfaces are normal congruences to some surfaces of
$\R^3$. We get as a corollary that a developable surface of
$\R^3$ is a critical point of the second-order functional 
$\mathcal{F}(S):=\int_S \sqrt{H^2-K} dA$. Things work exactly in the same way with 
$T \H^2,$ which can be identified with the set of oriented time-like lines of the 
Minkowski three-space $\R^{2,1}$, and whose Lagrangian surfaces are normal congruences
to space-like surfaces. We thus get that a developable space-like surface of 
 $\R^{2,1}$ is a critical point of the functional equivalent to ${\cal F}$
in $\R^{2,1}$.
In a forthcoming paper we shall
study more deeply the Hamiltonian stationary Lagrangian surfaces
of $(T\Sigma,\J,\G,\Om)$.

\bigskip

The paper is organised as follows: in Section 1 we give some
preliminary results and the precise statements of
 the two main theorems.
Section 2 is devoted to the proof of Theorem \ref{one}. The last
two sections deal with special cases: in Section 3 consider the
Euclidean case and prove Theorem \ref{two}; in Section 4 we take a
closer look to the special cases $T\S^2$ and $T\H^2$.

\bigskip

Finally the Authors wish to mention recent related results 
in the special Lagrangian case obtained independently by Dong \cite{Dg}. 

\section{Preliminaries and statements of results}

\subsection{The structures of $T\Sigma$}

In the following we consider an oriented Riemannian surface
$(\Sigma,g)$ and denote by $j$ the canonical complex structure
associated to it. We denote by $\pi$ the canonical projection of
the tangent bundle $T\Sigma$ onto its base $\pi: T \Sigma \to
\Sigma$. The two-dimensional subbundle $Ker (d\pi)$ of $TT \Sigma$
(it is thus a bundle over $T\Sigma$) will be called \em the
vertical bundle \em and denoted by $V\Sigma$.

We observe that we have not used the metric $g$ so far. The next
step consists of using the Levi--Civita connection $\nabla$ of $g$
to define the \emph{horizontal bundle} $H \Sigma$ as
follows: let $X$ be a tangent vector to $T\Sigma$ at some point
$(p_0,V_0)$. This implies that there exists a curve $\al(s)=
(p(s),V(s))$ such that $(\ga(0),V(0))=(p_0,V_0)$ and $\al'(0) = X$. If
$X \notin V \Sigma$ (which implies $p'(0) \neq 0$), we define the
connection map (\em cf. \em \cite{Ko},\cite{Do}) $K : TT \Sigma \to T\Sigma$ by $KX =
\nabla_{p'(0)}V(0),$ which does not depend on the curve $\al$.
If $X$ is vertical, we may assume that the
curve $\al$ stays in a fiber so that $V(s)$ is a curve in a
vector space. We then define $KX$ to be simply $V'(0)$. The
horizontal bundle is then $Ker (K)$ and we have a direct sum
\begin{equation}
\begin{array}{rcc} TT\Sigma= H\Sigma \oplus V\Sigma & \simeq & T\Sigma \oplus T\Sigma
 \\ X& \simeq &  (PX, KX )
  \end{array} \label{eq:directsum}
\end{equation}
Here and in the following, $P$ is a shorthand notation for $d\pi$.
We refer to \cite{Ko} and \cite{Do} for a more complete description of the horizontal and vertical bundles.
\bigskip

We shall use again the metric $g$ in order to pull back the canonical symplectic
form of $T^* \Sigma$ to a symplectic form $\Om$ in $T\Sigma,$ which admits a nice expression
in terms of the direct decomposition of $TT\Sigma$:
\begin{lemm} \label{symplectic}
Let $X$ and $Y$ be two tangent vectors to
$T \Sigma$; we have
 $$ \Om(X,Y) := g(K X, P Y)-g(PX,KY).$$
\end{lemm}
A proof of this lemma can be found in \cite{La}, p. 89.
\bigskip

We recall then the classical
\begin{defi}
A surface $L$ of $(T \Sigma,\Om)$ is said to be \em Lagrangian \em if
the restriction of $\Om$ vanishes on it.
\end{defi}

Next we define an almost complex structure $\J$ by $\J= j \oplus j,$ using the
direct sum $(\ref{eq:directsum})$ and the pseudo-metric $\G$ by the formula $ \G(.,.)=\Om(\J.,.)$. 
In \cite{GK1} it has been proved that $\G$ is a pseudo-Riemannian metric with signature $(2,2)$.
Proposition~\ref{Jcomplex} below shows that $\J$ is actually a complex
structure.

\subsection{Statements of the main theorems}

The projection map $\pi: T\Sigma \to \Sigma$ plays a crucial role in the local
classification of minimal 
Lagrangian surfaces of $(T\Sigma,\J,\G,\Om)$. If $L$ is some surface of $T\Sigma$ (not necessarily
Lagrangian), then the rank of the restriction to $L$ of the projection $\pi$ can be 0, 1 or 2 and is
locally constant. The case of rank 0 corresponds to the trivial
case of $L$ being a piece of a vertical fibre.

A simple example of Lagrangian surface of rank 1 is the 
\em
normal bundle \em over some curve $\ga$ of $\Sigma,$ i.e. the set of its normal lines to the curve $\ga$.
More precisely, denoting by $\vec{n}(s)$ a unit normal vector to
the curve at the point $\ga(s),$ the normal bundle of $\ga$ is the image of the immersion $X(s,t)=(\ga(s), t \vec{n}(s) )$. 
One can slightly generalize the construction by considering affine lines, i.e. adding
a translation term to the second factor of the immersion: $X(s,t)= (\ga(s), a(s) \vec{t} + t \vec{n}(s) ),$
where $\vec{t}$ denotes the unit tangent vector to $\ga(s)$ and $a(s)$ is some real-valued map. 
We shall call the image of such an immersion, which is still Lagrangian, 
an \em affine normal bundle over $\ga$. \em Affine normal bundles and their higher dimensional equivalents have been
introduced in the flat case in \cite{HL}, where they were called \em degenerate projections. \em

As the metric $\G$ is neutral, the induced metric on a surface of $(T\Sigma,\J,\G,\Om)$ may be degenerate. It is for example
the case of a vertical fibre and of the zero section $L_0:=\{ (p,0), p \in \Sigma \} \subset T\Sigma$. Such surfaces are
called \em null. \em

\medskip

 The first main result of this article characterizes rank one minimal Lagrangian surfaces and
 shows that, beyond the null surfaces, there is no rank two Lagrangian minimal surface if $(\Sigma,g)$ is non-flat:
 
\begin{theo} \label{one} Let $L$ be a smooth, non-null minimal Lagrangian
surface of $(T \Sigma,\J,\G,\Om)$. Then
\begin{itemize}
\item[(i)] either $L$ is the normal bundle over a geodesic on $(\Sigma, g)$, or 
\item[(ii)] $(\Sigma, g)$ is flat.
\end{itemize}
\end{theo}

The Lagrangian assumption in the theorem above is crucial: the existence of families of 
(non-Lagrangian) minimal surfaces in $(T \Sigma,\J,\G,\Om)$
has been proved in \cite{GK2}.
When the surface $\Sigma$ is flat, the situation appears to be richer, in the sense that
there exist many rank two minimal Lagrangian surfaces. As our classification is local,
there is no loss of generality to restrict ourselves to the Euclidean plane.

\begin{theo} \label{two} In the case where $(\Sigma,g,j)$
is the Euclidean plane $\R^2$ endowed with its
canonical inner product $\<.,.\>,$ the metric $\G$
on $T\R^2 \simeq \R^4$ is the flat pseudo-metric
of signature $(2,2)$. Moreover, if
 $L$ is a rank two minimal Lagrangian surface of $(T \R^2,\J,\G,\Om),$
 then 
 it is
parametrized by $X(p)=(p, \nabla u(p)),$ where the real map $u$ 
takes the following form
 $$ u(p)=f_1(\<p, V\>)+f_2(\<p,jV\>),$$
where $V$ is some constant unit vector of $\R^2$ and $f_1$ and $f_2$ are two
non-constant functions of the real variable of class $C^2$.
\end{theo}
In Section 3 we shall give, along with the proof of the theorem, a
geometric interpretation of the vector $V$.

\subsection{Some preliminary results}

We start with a result from \cite{Ko} which will be useful:
\begin{lemm} \cite{Ko} \label{bracket}
Given a vector field $X$ on $(\Sigma,g),$ there exists exactly one vector field $X^h$ and one vector field
$X^v$ on $T\Sigma$ 
(\footnote{The Reader should be aware that the notation for $X^v,X^h$ in Lemma \ref{bracket}
follows \cite{Ko} and corresponds to particular lifts of a field on $M$, whereas \cite{La}, in the proof of Lemma~\ref{symplectic}, uses the same notation to denote \emph{projections} of a vector field.}) 
such that $(PX^h,KX^h)=(X,0)$ and $(PX^v,KX^v)=(0,X)$.
Moreover, given two vector fields $X$ and $Y$ on $(\Sigma,g)$, we have, at the point $(p,V)$:
$$ [X^v,Y^v]=0$$
$$ [X^h,Y^v] \simeq (0,\nabla_X Y )$$
$$ [X^h,Y^h] \simeq ([X,Y], -R(X,Y)V),$$
where $R$ denotes the curvature of $g$ and we use the direct sum notation $(\ref{eq:directsum})$.
\end{lemm}

We say that a vector field $X$ on $T\Sigma$ is \em projectable \em if it is constant on the fibres. According to the 
lemma above, it is equivalent to the fact that there exists two vector fields
$X_1$ and $X_2$ on $\Sigma$ such that $X=(X_1)^h + (X_2)^v$.

\bigskip

We can now prove
\begin{prop} \label{Jcomplex}
The almost complex structure $\J$ is complex. 
\end{prop}
\emph{Proof.}
We compute the Nijenhuis tensor
\[
N(X,Y) = [X,Y] + \J [\J X,Y] + \J [X,\J Y] - [\J X,\J Y].
\]
Since $N(X,Y)$ depends pointwise on the tangent vectors we may assume for
computational purposes that $X$ and $Y$ are projectable, and use
lemma~\ref{bracket} together with the definition $\J=j \oplus j$. 
By linearity and skew-symmetry it suffices to prove that $N(X,Y)=0$ in
three distinct cases:
\begin{enumerate}
  \item vertical fields
  \[
	N(X^v,Y^v) = 0 
  \]
 \item horizontal fields
\begin{eqnarray*}
  N (X^h, Y^h) & = & - \Big( 0, R (X, Y) V + jR (jX, Y) V
  \\
  && + jR (X, jY) V - R(jX, jY) V \Big)
  \\
  & = & - \left( 0, jR (jX, Y) V + jR (X, jY) V \right)=(0,0)
\end{eqnarray*} 
 \item mixed fields
  \begin{eqnarray*}
	N (X^h, Y^v) & = & \left( 0, \nabla_X Y + j \nabla_{jX} Y + j \nabla_X (jY)
  - \nabla_{jX} (jY) \right)\\
  & = & \left( 0, \nabla_X Y + j \nabla_{jX} Y - \nabla_X Y - j \nabla_{jX} Y
  \right) = (0, 0)
  \end{eqnarray*}  
\end{enumerate}
where we have used the properties of K\"ahler manifolds: $\nabla j=0$, 
$R (j X, j Y) = R (X, Y)$.
\QED

\begin{coro}
The triple $(\G,\J,\Omega)$ defines a pseudo-K\"ahler structure on $T\Sigma$.
In particular $\J$ is parallel for the Levi-Civita connection.
\end{coro}

\bigskip

The following lemma describes the Levi-Civita connection $D$ of $\G$ in terms of the direct 
decomposition of $TT\Sigma$.

\begin{lemm} \label{connection}
Let $X$ and $Y$ two vector fields and assume that $Y$ is projectable, then at the point $(p,V)$ 
we have
\[
D_X Y = \left( \begin{matrix}
\nabla_{P X } P Y
\\
\nabla_{P X } K Y -{\textstyle {\frac{1}{2}}}
\Big( R(PX , P Y)V -jR(V,jP X)P Y-j R(V,j P Y)P X\Big)
\end{matrix} \right)
\]
where we have used column vector notation to indicate the components in the direct 
sum~(\ref{eq:directsum}).
\end{lemm}
\emph{Proof.} 
We use Lemma \ref{bracket} together with the Koszul formula:
\begin{eqnarray*}
2\G(D_X Y, Z) &=& X\G(Y,Z)+Y\G(X,Z)-Z\G(X,Y)+\G([X,Y],Z) \\
&& -\G([X,Z],Y)-\G([Y,Z],X) 
\end{eqnarray*}
where $X$, $Y$ and $Z$ are three vector fields on $T\Sigma$.
From the fact that $[X^v,Y^v]$ and $\G(X^v,Y^v)$ vanish we have:
\begin{eqnarray*}
2\G(D_{X^v} Y^v, Z^v) &=& X^v \G(Y^v,Z^v)+Y^v\G(X^v,Z^v)-Z^v\G(X^v,Y^v) \\
&& +\G([X^v,Y^v],Z^v)-\G([X^v,Z^v],Y^v)-\G([Y^v,Z^v],X^v)=0 .
\end{eqnarray*}
Moreover, taking into account that $\G(Y^v,Z^h)$ and similar quantities are constant on the fibres, 
we obtain
\begin{eqnarray*}
2\G(D_{X^v} Y^v, Z^h) &=& X^v \G(Y^v,Z^h)+Y^v\G(X^v,Z^h)-Z^h\G(X^v,Y^v) \\
&& + \G([X^v,Y^v],Z^h)-\G([X^v,Z^h],Y^v)-\G([Y^v,Z^h],X^v) \\
& = & -\G(-(\nabla_Z X)^v,Y^v)-\G(-(\nabla_Z Y)^v,X^v)=0.
\end{eqnarray*}
From these last two equations we deduce that $D_{X^v} Y^v$ vanishes.
Analogous computations show that $ D_{X^v} Y^h$ vanishes as well.
From Lemma \ref{bracket} and the formula $[X,Y]=D_X Y-D_Y X$, we deduce that
 $D_{X^h} Y^v \simeq (0, \nabla_X Y)$.

Finally, introducing $jW=-jR(X,Y)V-R(V,jY)X-R(V,jX)Y,$
we compute that
$$\G(D_{X^h} Y^h, Z^h)=\frac{1}{2}g(jW,Z)$$
and
$$\G(D_{X^h} Y^h, Z^v)= g(jZ,\nabla_X Y),$$
from which we deduce that
\begin{multline*}
D_{X^h} Y^h = (\nabla_X Y,\frac{1}{2}W)
\\
= \left(\nabla_X Y,\frac{1}{2} ( -R(X,Y)V+jR(V,jX)Y+jR(V,jY)X) \right) .
\end{multline*}

The conclusion of the proof follows easily. \QED
 
\bigskip

The fact that $\J$ is parallel with respect to $D$ implies the following
useful result about the extrinsic geometry of Lagrangian surfaces; 
this fact is known to hold in a positive K\"{a}hler manifold,
\emph{cf.}~\cite{Ch}.
\begin{lemm} \label{symmetric}
Let $L$ be a Lagrangian surface of $T\Sigma$ and $X,Y$ and $Z$ three vector fields tangent to $L$.
Then 
\[
h(X,Y,Z):= \Om(X,D_Y Z)=\G(\J X, D_Y Z)
\]
defines a tri-symmetric tensor called the tensor of extrinsic curvature. 
\end{lemm}
\emph{Proof.}
Let $\mathrm{II}$ denote the second fundamental form of the immersion:
\[
h(X,Y,Z)=\G(\J X, D_Y Z)=\G(\J X,\mathrm{II}(Y,Z)),
\]
which proves the tensorial nature of $h$ as well as the symmetry with respect to its last two
variables. 
\begin{eqnarray*}
h(X,Y,Z) &=& Y \G(\J X, Z) - \G(D_Y \J X, Z) = -\G(\J D_Y X,Z) = \G(\J Z, D_Y X) \\
&=& h(Z,Y,X)
\end{eqnarray*}
using the Lagrangian hypothesis on $L$. \QED

\bigskip

\section{Proof of Theorem \ref{one}}

The proof of Theorem \ref{one} will result from Propositions \ref{p1} and \ref{p2}, dealing with
Lagrangian surfaces of rank one and two, respectively.

\subsection{Rank one Lagrangian surfaces}

\begin{prop} \label{p1}
A rank one Lagrangian surface $L$ of $(T\Sigma,\J ,\G,\Om)$ 
 is an affine normal bundle over a curve $\ga$ of $\Sigma$. It is
 moreover H-minimal and the induced metric on $L$ is flat. 
Finally, $L$ is minimal if and only if the base curve $\ga$ is a geodesic
of $(\Sigma,g)$.
\end{prop}
\emph{Proof.}
A surface $L$ of $T\Sigma$ with rank 1 projection may be parametrized
locally by
$$ \begin{array}{lccc} X :
 & U
 &\to& T \Sigma \\
& (s,t)& \mapsto &  (\ga(s), V(s,t)),
  \end{array}$$
where $\ga(s)$ is a regular curve in $\Sigma$ and $V(s,t)$ some
tangent vector to $\Sigma$ at the point $\ga(s)$. 
Without loss of generality, we may assume that $\ga$ is parametrized by arclength, so that
 $\{\ga'(s), j \ga'(s) \}$ is an
orthonormal frame of $T\Sigma$ along the curve $\ga$.
 Writing $V=a \ga' + b j \ga'$ and using the
Fr\'enet equation $\nabla_{\ga'} \ga'= kj \ga'$, where $k$ denotes the curvature of $\ga,$
we compute the first derivatives of the
immersion (here and in the following, a letter in subscript denotes partial differentiation with 
respect to the corresponding variable). Using the direct sum notation:
\begin{eqnarray*}
X_s &=& (P X_s, K X_s) = (\gamma', \nabla_{\gamma'} \gamma') = (\ga', (a_s-kb)\ga'+ (b_s+ka) j\ga')
\\
X_t &=& (0, a_t \ga'+ b_t j \ga').
\end{eqnarray*}
If the immersion is Lagrangian, the following must vanish:
$$\Om(X_s,X_t)=-g(\ga',a_t \ga' + b_t j \ga')=-a_t.$$
 It follows that $a$ must be a function of $s$.
Then we see that for fixed $s$, the map $ t \mapsto (\ga(s),a(s)
\ga'(s)+b(s,t) j \ga'(s))$ parametrizes a line segment ruled by $j
\ga'(s)$ in $T_{\ga(s)} \Sigma$ which may be reparametrized by
\[
t \mapsto (\ga(s), a(s)\ga'(s) + t j \ga'(s))
\]
which we assume henceforth. We have thus proved the first part of Proposition \ref{p1}.

We compute easily that $X_s = (\ga', (a'-kt) \ga'+ ak j \ga')$ and $ X_t =(0, j\ga')$.
We also observe that the vector field (defined along the surface) $X_t$
depends only on the variable $s$, thus it can be extended to a
global vector field which is projectable. It follows that
we can use Lemma \ref{connection} in order to compute:
\[
D_{X_s} X_t = ( 0, \nabla_{\ga'} j \ga' ) =(0, -k\ga'), \qquad
D_{X_t} X_t = ( 0, 0).
\]
In view of Lemma \ref{symmetric}, the symmetric tensor $h(X,Y,Z)$ has four independent components. We calculate:
\[
h_{112}=\Om(X_s,D_{X_s} X_t)=
\Om( (\ga',(a'-kt) \ga'+ ak j \ga') , (0, -k \ga') ) = k
\]
$$h_{122}=\Om(X_s,D_{X_t} X_t)=0  \quad \quad h_{222}=\Om(X_t,D_{X_t} X_t)=0.$$ 
(As will become clear in a moment, we do not need the expression of $h_{111}$.) 

It remains to compute the induced metric, which is given in the coordinates $(s,t)$ by
$$ \left( \begin{array}{cc} - 2ak & -1 \\
  -1 &0\end{array}
\right).$$
We are now in a position to get the expression of the mean curvature vector:
\[
\G(2 \vec{H},\J X_s)=\frac{h_{111} G + h_{122}E-2 h_{112} F }{ E G - F^2} = - 2k
\]
and
\[
\G(2 \vec{H},\J X_t)=\frac{h_{112}G+h_{222}E-2h_{122}F}{EG-F^2}=0.
\]
It follows that
 $$ \vec{H}=k \J X_t = (0,k j\ga') =(0, \ga''(s)),$$
so that $L$ is minimal if and only if $k$ vanishes, namely $\ga$ is a geodesic.
Moreover, the determinant of the induced metric being $-1,$
we have the following formula:
$$ \mbox{div} \J \vec{H} = \mbox{div} (-k \pa_t)= 0$$
and hence $L$ is always Hamiltonian stationary.
Finally, denoting by $\bar{\nabla}$ and $\bar{R}$ the Levi-Civita connection and the curvature of the induced metric, an easy computation shows that
$\bar{\nabla}_{\pa_s} \pa_t$ and $\bar{\nabla}_{\pa_t} \pa_t$ vanish, so that 
$$ \bar{R}(\pa_t, \pa_s) \pa_t =\bar{\nabla}_{\pa_t} \bar{\nabla}_{\pa_s} \pa_t-\bar{\nabla}_{\pa_s}\bar{\nabla}_{\pa_t} \pa_t=0,$$
which implies the flatness of $L$.
\QED

\subsection{Rank two Lagrangian surfaces}

\begin{prop} \label{p2}
A rank two Lagrangian surface $L$ of $(T\Sigma,\J ,\G,\Om)$ 
 is the graph of the gradient of a real map $u$ on $(\Sigma,g)$:
$$ L :=\{ (p, \nabla u(p)), p \in \Sigma \} \subset T \Sigma. $$
Moreover, if $L$ is minimal then $g$ is flat.
\end{prop}
\emph{Proof.}
A rank 2 surface is nothing but the graph of a vector field $V(p)$ of $\Sigma$ and thus is the image 
of the immersion $X(p)=(p,V(p))$. Let $(s,t)$ be conformal local coordinates on $(\Sigma, g)$ 
such that $j \pa_s = \pa_t$ and $j\pa_t = -\pa_s$. We denote by $r(s,t)$ the logarithmic
conformal factor, so that the metric takes the following form:
$g(s,t)= e^{2r}(ds^2+dt^2)$.
A standard computation shows that
\begin{eqnarray*}
\nabla_{\pa_s} \pa_s &=& r_s \pa_s - r_t \pa_t \\
\nabla_{\pa_t} \pa_s &=& \nabla_{\pa_s} \pa_t = r_t \pa_s + r_s \pa_t \\
\nabla_{\pa_t} \pa_t &=& -r_s \pa_s + r_t \pa_t \, .
\end{eqnarray*}
The following relations between the curvature tensor $R$, the Gauss curvature $K$ and the conformal factor $r$ 
of $(\Sigma,g)$ will be useful later:
$$ K=e^{-4r}g(R(\pa_s,\pa_t)\pa_s,\pa_t)=-e^{-2r}\Delta r.$$
 Writing $V(s,t)=P(s,t)\pa_s + Q(s,t)\pa_t$, the first derivatives of the immersion are:
$$ X_s=(\pa_s, (P_s + Pr_s +Qr_t)\pa_s + (Q_s- Pr_t + Qr_s)\pa_t),$$
$$ X_t=(\pa_t, (P_t + Pr_t -Qr_s)\pa_s + (Q_t+ Pr_s + Qr_t)\pa_t),$$
so that
\begin{eqnarray*}
\Om(X_s,X_t) &=& g((Q_s- Pr_t + Qr_s)\pa_t,\pa_t)-g(\pa_s,(P_t + Pr_t -Qr_s)\pa_s) \\
&=& e^{2r}(Q_s+2Qr_s - P_t - 2Pr_t).
\end{eqnarray*}
Thus the Lagrangian condition is equivalent to $ (Pe^{2r})_t=(Qe^{2r})_s,$ so that there exists locally
a real map $u$ on $\Sigma$ such that $Pe^{2r}=u_s$ and $Qe^{2r}=u_t$; in other words, the vector field
$V$ is the gradient of $u,$ and we have the first part of Proposition \ref{p2}.

\bigskip

Next a parametrization of $L$ is
$$ \begin{array}{lccc} X :
 & \Sigma
 &\to& T \Sigma \\
& (s,t)& \mapsto &  (p(s,t), e^{-2r}(u_s \pa_s+u_t \pa_t)),
  \end{array}$$
and we compute
\begin{eqnarray*}
X_s &=& (\pa_s, \nabla_{\pa_s} \nabla u) \\
&=& (\pa_s, e^{-2r} \left( (u_{ss}-2 r_s u_s)\pa_s + u_s \nabla_{\pa_s} \pa_s 
+ (u_{st}-2 r_s u_t)\pa_t + u_t \nabla_{\pa_s} \pa_t) \right) \\
&=& (\pa_s, e^{-2r}(u_{ss}- r_s u_s+r_t u_t)\pa_s + e^{-2r}(u_{st}-r_s u_t - r_t u_s)\pa_t )
\end{eqnarray*}
Analogously
\[
X_t =(\pa_t, e^{-2r}(u_{st}- r_s u_t-r_t u_s)\pa_s + e^{-2r}(u_{tt}-r_t u_t + r_s u_s)\pa_t ) .
\]
 Denoting for simplicity
 $$X_s:=(\pa_s, a \pa_s + b \pa_t)  \quad \quad X_t:=(\pa_t, b \pa_s + c \pa_t ),$$
the induced metric is given by
\[
E = \G(X_s,X_s)=\Om(\J X_s, X_s)=g(j( a \pa_s + b \pa_t),\pa_s ) - g(j\pa_s, a \pa_s + b \pa_t ) = -2 b e^{2r}, 
\]\[
F=\Om(\J X_s, X_t)=g(j (a \pa_s + b \pa_t),\pa_t )- g(j\pa_s, b \pa_s + c \pa_t )= (a-c) e^{2r},
\]\[
G = \Om(\J X_t, X_t)=g(j ( b \pa_s + c \pa_t),\pa_t ) - g(j\pa_t, b \pa_s + c\pa_t )= 2 b e^{2r}.
\]

Moreover, the vector fields $X_s$ and $X_t$ admit extensions on $T\Sigma$ which
are projectable, to that we can use Lemma \ref{connection}, to get
\[
D_{X_s}X_s =\Big( r_s \pa_s - r_t \pa_t,
(a_s +ar_s +br_t) \pa_s+(b_s- ar_t+br_s)\pa_t+u_s e^{-2r} j R(\pa_s,\pa_t )\pa_s \Big)
\]
\begin{multline*}
D_{X_{s}} X_t =
\Big(r_t \pa_s + r_s \pa_t, (b_s+ b r_s + c r_t) \pa_s+ (c_s -b r_t + c r_s) \pa_t\Big) \\
+\frac{1}{2} \Big( 0,-R(\pa_s,\pa_t)\nabla u + u_s e^{-2r} j R(\pa_s,\pa_t )\pa_t 
+ u_t e^{-2r} j R(\pa_s,\pa_t)\pa_s \Big),
\end{multline*}
\begin{multline*}
D_{X_t} X_t = 
\Big( -r_s \pa_s + r_t \pa_t, (b_t + b r_t -c r_s) \pa_s + (c_t + b r_s + c _t) \pa_t
- u_t e^{-2 r} j R(\pa_t, \pa_s) \pa_t \Big) \, .
\end{multline*}
This allows us to calculate the following components of the tensor $h$:
\begin{eqnarray*}
h_{111} &=& h(X_s,X_s,X_s)=\Om(X_s,D_{X_{s}} X_s) \\
&=& g(a\pa_s + b \pa_t,r_s \pa_s - r_t \pa_t) \\
&& -g(\pa_s,(a_s +ar_s +br_t) \pa_s) - u_s e^{-2 r} g(\pa_s, jR( \pa_s ,\pa_t )\pa_s) \\
&=&  e^{2r} (ar_s - br_t -(a_s+ar_s +br_t) ) + u_s e^{-2r}g(\pa_t,R(\pa_s,\pa_t )\pa_s) \\
&=& e^{2r} (-a_s-2br_t+u_s K).
\end{eqnarray*}
and similarly\footnote{
Note that coefficients $h_{112}$ and $h_{122}$ can be computed by two different
methods, yielding two seemingly different expressions.} with the other
coefficients of $h$. Consequently
\begin{multline*}
\G(2\vec{H},\J X_s)=\frac{h_{111}G + h_{122}E - 2 h_{112}F}{EG-F^2}
=\frac{2b(h_{111}-h_{122}) - 2(a-c)h_{112}}{-e^{2r}(4b^2 +
(a-c)^2)} 
\\
=\frac{2b((a-c)_s + 4br_t) + 2(a-c)(-b_s+ (a-c) r_t) )}{4b^2 +(a-c)^2}
\\
= \frac{(a-c)_s (2b)-(a-c)2b_s}{(2b)^2 + (a-c)^2 } +2r_t 
= (\arg(2b+i (a-c))_s +2r_t \, .
\end{multline*}

A similar computation yields
$$\G(2\vec{H},\J X_t)=(\arg(a-c)+2ib))_t - 2r_s,$$
and hence the vanishing of $\vec{H}$ implies
\begin{eqnarray*}
(\arg(c-a +2ib))_s - 2r_t=0
\\
(\arg(c-a+2ib))_t + 2r_s=0 .
\end{eqnarray*}
Differentiating the first equation with respect to the variable $t$, and 
the second equation with respect to the variable $s$ yields $\Delta r = 0$,
which implies that $\Sigma$ has vanishing curvature and concludes the proof 
of Proposition \ref{p2}. \QED

\section{Minimal Lagrangian surfaces in $T \R^2$}

 We now consider the Euclidean plane $(\R^2,\<.,.\>)$ with coordinates $(x_1,x_2);$
the metric is $dx_1^2+dx_2^2$ and
 the complex structure is $j(x_1,x_2)=(-x_2,x_1)$.
On $T \R^2$ we define the coordinates
$(x_1,x_2,y_1,y_2), $ in which the canonical symplectic structure writes
$$ \Om = dy_1 \wedge dx_1 + dy_2 \wedge dx_2.$$
The complex structure is $\J := j \oplus j$ and the metric $\G$ by
$$ \G= dx_2 dy_1 - dx_1 dy_2.$$

\bigskip

The complex structure $\J $ induces an identification of $T\R^2 \simeq \R^{2,2}$ with $\C^2$, given
by $(x_1,x_2,y_1,y_2) \simeq (w_1:=x_1 + i x_2, w_2:=y_1+iy_2)$.
The pseudo-Hermitian metric takes the form:
$$ H(.,.)=\G(.,.) + {i} \Om(.,.)=\G(.,.) + {i} \G(.,\J .).$$
 
 However we can consider on $\C^2\simeq T\R^2$ the canonical Riemannian and symplectic 
 structures, and define classically the Lagrangian angle: if $e_1 \wedge e_2$ denotes 
 a Lagrangian plane, its Lagrangian angle $\beta$ is the argument of
 $dw_1 \wedge dw_2 (e_1 \wedge e_2)=\det_{\C} (e_1,e_2)$. If $L$
 is some Lagrangian surface, the Lagrangian angle function $\beta$ is defined on $L$ 
 by $\beta(p)=\beta(T_p L)$. The Reader should note that $\beta$ bears a priori no relation to 
 the pseudo-K\"{a}hlerian structure defined on $T\R^2$. Nevertheless one obtains the following surprising 
 
\begin{prop}
The relation
\[
\vec{H} =\frac{1}{2} \J D \beta
\]
still holds for Lagrangian surfaces of $\R^{2,2}$, where $D\beta$ denotes the gradient of $\beta$ in the induced (pseudo-)metric. 

In particular, a Lagrangian surface is minimal
if and only if its Lagrangian angle is (locally) constant.
\end{prop}
\emph{Proof}. Let $(e_1,e_2)$ a frame along $L$ such that 
$\G(e_1,e_1)=-1,\G(e_2,e_2)=1$ and $\G(e_1,e_2)=0$. The Lagrangian assumption implies $H(e_1,e_1)=-1$, $H(e_2,e_2)=1$ and $H(e_1,e_2)=0$. Thus, given a vector 
$\vec{V}$ of $\R^4,$ the following formula holds:
\[
\vec{V} = -H(\vec{V},e_1) e_1 + H(\vec{V},e_2)e_2.
\]
We differentiate the relation $e^{i \beta(p)} = \det_{\C}(e_1(p),e_2(p))$ 
with respect to $e_1$, which yields, using the fact that $dw_1 \wedge dw_2$ is parallel 
with respect to the Levi-Civita connection induced by $\G$:
\begin{eqnarray*}
i e_1( \beta) e^{i \be} &=& 
 \det_{\C}(D_{e_1} e_1,e_2)+\det_{\C}(e_1,D_{e_1} e_2) 
 \\
&=& - H(D_{e_1} e_1,e_1)\det_{\C}( e_1,e_2) + H(D_{e_1} e_2,e_2)\det_{\C}( e_1,e_2)
\\
&=& e^{i \beta} \big[ (-\G(D_{e_1} e_1,e_1)+\G(D_{e_1} e_2,e_2))
\\
&& \qquad + i (-\G(D_{e_1} e_1,\J e_1)+\G(D_{e_1} e_2,\J e_2)) \big] \, .
\end{eqnarray*}
Thus $e_1(\beta) = - h(e_1,e_1,e_1) + h(e_1,e_2,e_2)= \G(2 \vec{H}, \J e_1)$, 
proving that 
\[
\G(\J D \be, \J e_1)=\G(2\vec{H},\J e_1) .
\]
Analogously we prove that $ \G(\J D \be, \J e_2)=\G(2\vec{H},\J e_2)$ 
and the proof is complete. \QED

\begin{rema}
The surprising fact that one uses the same definition for the Lagrangian angle 
though the underlying (pseudo-)K\"ahler is quite different can be explained by looking 
at the isometries for that structure. Indeed the group is $(\G,\J ,\Om)$-preserving matrices is 
none other than $U(1) \times SL(2,\R)$, written in complex notations as $2 \times 2$ matrices.
(While the corresponding group in flat $\C^2$ is $U(2) = U(1) \times SU(2)$.) So that 
in both cases the Lagrangian angle measures the $U(1)$ factor.
\end{rema}

 \bigskip

 We are now in position to determine locally the minimal Lagrangian surfaces of $\R^{2,2}$:

\begin{prop}
Let $L$ be a rank two Lagrangian surface of $\R^{2,2},$ i.e. it is
parametrized by $X(p)=(p, \nabla u(p)),$ where $u$ is a $C^2$
map defined on an open subset of $(\R^2,\<.,.\>).$ Then $L$ has
constant Lagrangian angle $\beta_0$ if and only if it takes the
following form
 $$ u(p)=f_1(\<p, e^{i (\beta_0/2+\pi/4)}\>)+f_2(\< p, ie^{i (\beta_0/2+\pi/4)}\>),$$
where $f_1$ and $f_2$ are two non-constant functions of the real
variable of class $C^2$.
\end{prop}
\emph{Proof.}
We first compute the first derivatives of the immersion $X,$ writing $p=(s,t)$:
$$ X_s=(1,0,u_{ss},u_{st}) \simeq (1, u_{ss} +iu_{st})
\hspace{2em} X_t=(0,1,u_{st},u_{tt}) \simeq (i,u_{st}+ iu_{tt}),$$
so the Lagrangian angle map is given by:
$$ \beta(s,t) = \arctan \left( \frac{u_{tt}- u_{ss}}{2 u_{st}} \right)$$
and the constant Lagrangian angle condition translates into
the linear PDE
$$\cos \beta_0(u_{tt}- u_{ss}) - 2 \sin \beta_0 u_{st} =0,$$
In order to solve this, we introduce the linear change of variables defined by
$$\left( \begin{array}{c} \si \\ \tau \end{array} \right) =
 \left( \begin{array}{cc} \cos \te & \sin \te \\
  -\sin \te & \cos \te \end{array}
\right) \left( \begin{array}{c} s \\ t \end{array} \right)
,$$
where $\te$ is some fixed constant, 
so that
$$ u_{tt}=\sin^2 \te u_{\si \si} + \cos^2 \te u_{\tau \tau} + 2 \cos \te \sin \te u_{\si \tau},$$
$$u_{ss}=\cos^2 \te u_{\si \si} + \sin^2 \te u_{\tau \tau} - 2 \cos \te \sin \te u_{\si \tau},$$
$$ u_{st}= (\cos^2 \te - \sin^2 \te)u_{\si \tau}+ \cos \te \sin \te (u_{\si \si}-u_{\tau \tau}),$$
and thus
\begin{multline*}
\cos \beta_0(u_{tt}- u_{ss}) - 2 \sin \beta_0 u_{st} =
\left(\cos \beta_0(\sin^2 \te-\cos^2 \te) - 2 \sin \beta_0 \cos \te \sin \te \right)u_{\si \si} 
\\
+\left(\cos \beta_0(\cos^2 \te -\sin^2 \te)+  2 \sin \beta_0 \cos \te \sin \te \right)u_{\tau \tau}
\\
+ \left( 4\cos \beta_0 \cos \te \sin \te - 2 \sin \beta_0 (\cos^2 \te -\sin^2 \te) \right) u_{\si \tau} 
\\
 = \cos(2\te-\be_0)(u_{\tau \tau}-u_{\si \si}) + 2 \sin (2\te-\be_0)u_{\si \tau}.
\end{multline*}
Hence, choosing $\te=\be_0/2 + \pi/4,$ the equation becomes $u_{\si \tau}=0$,
whose general solution is
\begin{eqnarray*}
u(s,t) &=& f_1(\si)+f_2(\tau)=f_1(\cos \te s+ \sin \te t)+f_2(-\sin \te s+ \cos \te t)
\\
& =& f_1(\< p, e^{i \te}\>)+f_2(\< p, ie^{i \te}\>).
 \end{eqnarray*} \QED
\bigskip
It is easy to see that the second part of Theorem \ref{two} is essentially a rewriting of the previous proposition.

\begin{rema} The formula $ \beta = \arctan \left( \frac{u_{tt}- u_{ss}}{2 u_{st}} \right)$
might be compared with the one we have in the classical
(Riemannian) case:
$$ \beta = \arctan \left( \frac{\Delta u}{1 - \det Hess(u)} \right).$$
In the classical case, the minimality is expressed by a kind of ``interpolation" 
between the Laplace and Monge-Amp\`ere equation.
Here, we can regard the equation as an interpolation between two
hyperbolic equations, the wave equation and the operator
$\pa_{st}$.
\end{rema}

\begin{rema}
The fact that here minimal Lagrangian surfaces can be only of class of $C^1$ makes a
great contrast with the positive case, where Special Lagrangian
 surfaces must be analytic (the underlying equation being elliptic).
\end{rema}

\begin{exem} Taking for example $u(s,t)=\sin s + \cos t,$ we get a doubly
periodic minimal Lagrangian surface in $\R^4$ or equivalently a
compact minimal Lagrangian surface in $T \mathbb{T}^2$.
\end{exem}

\section{The case of $T \S^2$ and normal congruences of surfaces in $\R^3$}

It is well known that the normal congruence to a regular, oriented surface $S$ of $(\R^3,\<.,.\>)$
defines a Lagrangian surface $\bar{S}$ in the space $\mathbb L^3$ of oriented lines of $\R^3$.
The latter is naturally identified with $T\S^2$ by the following
\[
\mathbb L^3 \ni \{ V+tp, t \in \R \} \simeq (p,V-\<V,p\>p) \in T \S^2 .
\]
Since $T\S^2$ is naturally and isometrically embedded in $T\R^3 = \R^3 \times \R^3$ as the submanifold
\[
\mathcal{S} = \{ (N,Y) \in \R^3 \times \R^3 , \, \<N,Y\>=0 \},
\]
we have two ways of describing a tangent vector $\xi$ at a point $(N,Y)$:
\begin{itemize}
\item it can be seen as a couple $\xi \simeq (\nu,\eta)$ in $\R^3 \times \R^3$ such that 
$\< N,\nu \> = 0$ and $\< N, \eta \> + \< \nu, Y \> = 0$, or
\item using the direct sum (\ref{eq:directsum}), we see that $P \xi = \nu$ and 
$K \xi = \nabla_\nu Y$ where $s \mapsto Y(s)$ extends along the tangent direction $\eta$; then 
$K \xi$ is the tangential projection $(\eta)^T = \eta - \< \eta, N \> N$.  
\end{itemize}

\begin{lemm} \label{h-mini}
The deformation of a regular surface $S$ of $\R^3$ induces a Hamiltonian deformation of $\bar{S}$ in
$T\S^2$.
\end{lemm}
\emph{Proof.} Let $X: U \to \R^3$ a local
parametrization of $S,$ $N$ the unit normal vector field and $h$ a
compactly supported function on $U$; we furthermore assume that
$X$ is a parametrization along the lines of curvatures, so that,
denoting by $\la$ and $\mu$ the curvature functions, we have the two
equations
$N_s=\la X_s$ and $N_t= \mu X_t$.

We consider a normal variation $V=hN,$ where $h$ is some smooth real map on $U$. Starting from
$X^\eps=X+ \eps h N,$ we have
$$ X_s^\eps=X_s + \eps (h_s N + h N_s) \hspace{2em} X_t^\eps=X_t + \eps (h_tN+hN_t),$$
so that
$$ X_s^\eps \times X_t^\eps =X_s \times X_t + \eps W+ o(\eps),$$
where
\[
W:= h_s N \times X_t+ h N_s\times X_t + h_t X_s \times N + h X_s \times N_t .
\]
Consequently
\[
|X_s^\eps \times X_t^\eps | = |X_s \times X_t | + \eps \<N ,W\> + o(\eps)
 \]
and
\[
N^\eps = N + \eps \, \frac{W - \< W,N \> N}{|X_s \times X_t |}+o(\eps) 
= N + \eps \, \frac{W^T}{|X_s \times X_t |}+o(\eps),
\]
where $W^T$ denotes again the tangential projection.

Introducing the notations $e_1:=X_s/|X_s|$ and $e_2:=X_t/|X_t|$, we have
\[
W^T= h_s N \times X_t+ h_t X_s \times N = h_s |X_t| N \times e_2+h_t|X_s|e_1 \times N ,
\]
so that
\[
N^\eps = N - \eps \left( \frac{h_s}{ |X_s|} e_1 + \frac{h_t}{|X_t|} e_2 \right) + o(\eps).
\]

The next step consists of looking at the effect of this normal variation \hbox{$V=hN$} on the normal congruence.
A parametrization of the normal congruence of $X$ being $\bar{X}=(N,X-\<X,N\>N),$ we have
\begin{multline*}
X^\eps - \< X^\eps,N^\eps\>N^\eps = X+ \eps hN 
\\
 - \left\< X+\eps hN, N - \eps \left( \frac{h_s}{ |X_s|} e_1 + \frac{h_t}{|X_t|} e_2 \right) \right\>
 \left(N - \eps \left( \frac{h_s}{ |X_s|} e_1 + \frac{h_t}{|X_t|} e_2 \right) \right) + o(\eps)
 \\
 = X - \<X,N\>N + \eps \left[ \left\< X, \frac{h_s}{ |X_s|} e_1 + \frac{h_t}{|X_t|} e_2 \right\> N
 + \< X,N \> \left( \frac{h_s}{ |X_s|} e_1 + \frac{h_t}{|X_t|} e_2 \right) \right] 
 \\
 + o(\eps)
\end{multline*}
so finally $\bar{X}^\eps=\bar{X} + \eps \bar{V} + o(\eps)$ with
\begin{multline*}
\bar{V} = \Bigg( -\frac{h_s}{ |X_s|} e_1-\frac{h_t}{|X_t|} e_2,
\\
\left\< X, \frac{h_s}{ |X_s|} e_1 + \frac{h_t}{|X_t|} e_2 \right\> N
 + \< X,N \> \left( \frac{h_s}{ |X_s|} e_1 + \frac{h_t}{|X_t|} e_2 \right) \Bigg)
\end{multline*}
where we have used the $T\mathcal{S}$ formalism. So that
\[
P \bar{V} = -\frac{h_s}{ |X_s|} e_1-\frac{h_t}{|X_t|} e_2
\textrm{ and }
K \bar{V} = \< X,N \> \left( \frac{h_s}{ |X_s|} e_1 + \frac{h_t}{|X_t|} e_2 \right) \; .
\]
In order to understand the \em normal \em variation induced by $\bar{V}$ on $\bar{S},$ we compute
a basis of its normal space.
\[
P \bar{X}_s = N_s = \la |X_s|e_1 
\]
\begin{eqnarray*}
K \bar{X}_s &=& X_s - \big( \<X_s,N\>N-\<X,N_s\>N-\<X,N\>N_s \big)^{T} 
\\
&=& X_s-\<X,N\>N_s = (1-\la \<X,N\>) |X_s| e_1
\end{eqnarray*}
Analogously, we have
\[
P \bar{X}_t = \mu |X_t|e_2 \; , \quad K \bar{X}_t (1-\mu \<X,N\>)|X_t|e_2 
\]
It is then obvious to compute the orthonormal basis for the normal bundle and we deduce
\[
\G(\bar{V},\J \bar{X}_s)= \Om(\bar{V},\bar{X}_s) = \<X,N\>\la h_s + h_s (1-\la \<X,N\>)=h_s
\]
and similarly $\G(\bar{V},\J \bar{X}_t)=h_t$. This means that $ \bar{V}^{\perp}=\J D h$, 
i.e. the vector field $ \bar{V}^{\perp}$ is Hamiltonian. \QED

\begin{lemm} \label{corresp} Let $S$ be a surface in $\R^3$ and $\bar{S}$ its
normal congruence. We denote by $\mathcal{A}(\bar{S})$ the area with respect 
to the metric $\G$ and by $\mathcal{F}(S)$ the functional
 defined by $\mathcal{F}(S):=\int_S \sqrt{H^2-K} dA$,
where $H$ and $K$ are respectively the mean curvature and the Gauss curvature of
$S$. Then
\[
\mathcal{A}(\bar{S})=\mathcal{F}(S).
\]
\end{lemm}
\emph{Proof.}
From the expressions for $\bar{X}_s$ and $\bar{X}_t$ computed in the proof of Lemma \ref{h-mini},
we obtain the coefficients of the first fundamental form of the immersion~$\bar{X}$:
\[
\bar{E}=\bar{G}=0, \qquad
\bar{F}=\Om(\J \bar{X}_s,\bar{X}_t)=(\mu - \lambda)|X_s| |X_t|.
\]
It follows that
\begin{eqnarray*}
\int_U \sqrt{|\bar{E}\bar{G}-\bar{F}^2|}dsdt &=& 
\int_U |\bar{F}| dsdt = \int_U |\la -\mu| \sqrt{EG-F^2}dsdt
\\
&=& \int_{X(U)} \sqrt{H^2-K} dA
\end{eqnarray*}
so ${\cal A}(\bar{S})={\cal F}(S)$. \QED

\bigskip

Lemmas \ref{h-mini} and \ref{corresp}  prove that the normal congruence
of a surface $S$ of $\R^3$ is Hamiltonian
stationary if and only if $S$ is a critical point of ${\cal F}$. On the other hand,
we know by Proposition \ref{p1} that rank one Lagrangian surfaces of $T\S^2$ are Hamiltonian stationary.
The next lemma provides a geometric interpretation of the rank one condition:

\begin{lemm} A non-planar surface $S$ of $\R^3$ is developable if and only if
its normal congruence defines a rank one Lagrangian surface in
$T\S^2$, i.e. is the normal
bundle of some curve of $\S^2$.
\end{lemm}
\emph{Proof.} By definition, a developable surface has
vanishing Gauss curvature, which implies that the Gauss image is a
curve (or a single point) in $\S^2$. 
As the Gauss map of $S$ is
nothing but the projection of $\bar{S}$ on the base $\S^2,$ the
result follows.

\bigskip

Finally, putting all these facts together we get:

\begin{coro}
A developable surface of $S$ of $\R^3$ is a critical point of the functional ${\cal F}$.
\end{coro}

Finally everything in Section 4 can be readily adapted to the case of
$T \H^2$, which is identified with the set of positive time lines $\mathbb L_+^3$ 
of the Minkowski space $(\R^{2,1}, \<.,.\>_1),$ by the following
$$\mathbb L_+^3 \ni \{ V+tp, t \in \R \} \simeq (p,V-\<V,p\>_1 p) \in T \H^2.$$
Here, $\H^2$ denotes the hyperboloid model of the hyperbolic plane, that is the space-like quadric
$$\H^2:=\{ p \in \R^{2,1}, \<p,p\>_1=-1, p_3 >0 \}.$$
We leave to the Reader the easy task to check 
that a developable space-like surface of $(\R^{2,1} ,\<.,.\>_1)$ is a critical point of the functional equivalent to ${\cal F}$
in $\R^{2,1}$.

\bigskip \bigskip

\noindent
Henri Anciaux \\
Department of Mathematics and Computing \\
Institute of Technology, Tralee \\
Co. Kerry, Ireland \\
henri.anciaux@ittralee.ie \\

\medskip
\noindent
Brendan Guilfoyle \\
Department of Mathematics and Computing \\
Institute of Technology, Tralee \\
Co. Kerry, Ireland \\
brendan.guilfoyle@ittralee.ie \\

\medskip
\noindent
Pascal Romon \\
Universit\'{e} de Paris-Est Marne-la-Vall\'{e}e \\
5, bd Descartes, Champs-sur-Marne \\
77454 Marne-la-Vall\'{e}e cedex 2, France \\
pascal.romon@univ-mlv.fr

\end{document}